\documentclass[11pt,a4paper]{article}
\usepackage{amssymb,amsmath,color,graphicx}
\catcode`!=11
\let\!int\int \def\int{\displaystyle\!int}
\let\!lim\lim \def\lim{\displaystyle\!lim}
\let\!cap\cap \def\cap{\displaystyle\!cap}
\let\!cup\cup \def\cup{\displaystyle\!cup}
\let\!sup\sup \def\sup{\displaystyle\!sup}
\let\!inf\inf \def\inf{\displaystyle\!inf}
\let\!sum\sum \def\sum{\displaystyle\!sum}
\let\!max\max \def\max{\displaystyle\!max}
\let\!min\min \def\min{\displaystyle\!min}
\catcode`!=12
\let\oldsection\section
\renewcommand\section{\setcounter{equation}{0}\oldsection}

\newtheorem{lemma}{Lemma}[section]
\newtheorem{theorem}{Theorem}[section]
\newtheorem{remark}{Remark}[section]
\newtheorem{proposition}{Proposition}[section]

\newcommand{\R}{{\mathbb{R}}}

\allowdisplaybreaks

\def\d{\delta}

\begin{document}

\title{Global stability of traveling waves with oscillations for Nicholson's blowflies equation
\footnotetext{\small E-mails:  mei@champlaincollege.qc.ca,
 mei@math.mcgill.ca(M. Mei), \\ zhangkj201@nenu.edu.cn(K. Zhang), zhangqifeng0504@163.com(Q. Zhang)}}

\author{
Ming Mei$^{a,b}$,  Kaijun Zhang$^c$, and Qifeng Zhang$^d$
\\
\\
{ \small \it $^a$Department of Mathematics, Champlain College Saint-Lambert}
\\
{ \small \it Quebec,  J4P 3P2, Canada, and}
\\
{ \small \it $^b$Department of Mathematics and Statistics, McGill University}
\\
{ \small \it Montreal, Quebec,   H3A 2K6, Canada}
\\
{ \small \it $^c$School of Mathematics and Statistics, Northeast Normal University}
\\
{ \small \it Changchun, Jilin,   130024, China}
\\
{ \small \it $^d$School of Science, Zhejiang Sci-Tech  University}
\\
{ \small \it Hangzhou, Zhejian,   310018, China}
}

\maketitle

\medskip

\begin{abstract}
For Nicholson's blowflies equation, a kind of  reaction-diffusion equations with time-delay, when the ratio of birth rate coefficient and death rate coefficient satisfies $\frac{p}{\delta}>e$,  the large time-delay $r>0$ usually causes the traveling waves to be oscillatory. In this paper, we are interested in the global stability of these oscillatory traveling waves, in particular, the challenging case of the critical traveling waves with oscillations.   We   prove  that,  the  critical oscillatory traveling waves  are globally stable with the algebraic convergence rate  $t^{-1/2}$, and the non-critical traveling waves are globally stable with the exponential convergence rate  $t^{-1/2}e^{-\mu t}$ for a positive constant $\mu$,  where the initial perturbations around the oscillatory traveling wave in a weighted Sobolev can be  arbitrarily large. The approach adopted is the technical weighted energy method  with some new development in establishing  the boundedness estimate of the oscillating solutions, which, with the help of optimal decay estimates by deriving the fundamental solutions for the linearized equations, can allow us to prove the global stability and to obtain the optimal convergence rates.


 \end{abstract}

\pagestyle{myheadings} \thispagestyle{plain} \markboth{M. Mei, K. Zhang and Q. Zhang}{Global stability of oscillatory traveling waves}

\newpage

\section{Introduction and main result}\label{int}
This is a continuation of the previous studies \cite{CMYZ,LLLM} on the stability of oscillatory traveling waves for Nicholson's blowflies equation, a class of non-monotone  reaction-diffusion equations with time-delay, which describes the population dynamics of a single species like  the Australian blowflies
\cite{GSW, GW,Mei-Lin-Lin-So1,Mei-Lin-Lin-So2,MSLS,SY}:
\begin{equation}
\begin{cases}
\dfrac{\partial v(t,x)}{\partial t} - D\dfrac{\partial^2
v(t,x)}{\partial x^2} +d(v(t,x)) = b(v(t-r,x)), \ \ (t,x)\in
\R_+\times \R, \\
v(s,x)=v_0(s,x), \ \ s\in [-r,0], \ x\in \R.
\end{cases}
\label{1.1}
\end{equation}
Here, $v(t,x)$ is the
mature population at time $t$ and location $x$, $D>0$ the spatial diffusion rate of the mature blowflies.
$d(v)$ and $b(v)$ represent  Nicholson's death rate function and Nicholson's birth rate function, respectively, in the forms of
\begin{equation}
d(v)=\delta v, \ \ \ b(v)=pve^{-av},
\label{1.1-1}
\end{equation}
where $\delta>0$ is the death rate coefficient, $p>0$ is the maximal egg daily production rate per blowfly, $a>0$ is a constant, $r>0$ is the matured age of blowflies, the so-called time-delay.

Clearly, the equation \eqref{1.1} possesses two constant equilibria
\[
v_-=0 \mbox{ and }  v_+=\frac{1}{a}\ln \frac{p}{\delta}.
\]
When $\frac{p}{\delta}>1$, then $v_+>v_-=0$.
Throughout this paper, naturally we assume that
\begin{equation}
\lim_{x\to\pm\infty} v_0(s,x)=v_\pm \ \ \mbox{ uniformly in } s\in[-r,0].
\label{2017-1}
\end{equation}

A traveling wave for \eqref{1.1} is a
 special solution to \eqref{1.1}
of the form $\phi(x+ct)\ge 0$ with $\phi(\pm \infty) = v_{\pm}$:
\begin{equation}
\begin{cases}
c\phi'(\xi)-D\phi''(\xi) + \d\phi(\xi)=b(\phi(\xi-cr)), \\
\phi(\pm\infty)=v_\pm,
\end{cases} \label{1.3}
\end{equation}
where $\xi=x+ct$, $'=\frac{d}{d\xi}$, and $c$ is the wave speed.

The main purpose of the paper is to prove the global stability of the oscillatory traveling waves as well as the optimal convergence rates.

First of all, let us review the progress in the existence of traveling waves. When $1<\frac{p}{\delta}\le e$, by using the upper-lower solutions method, So and Zou \cite{SZ} first proved that there exists a minimal wave speed $c_*=c_*(r)>0$,
the so-called the critical wave speed, such that when $c\ge c_*$, the equation \eqref{1.3} possesses monotone traveling waves $\phi(x+ct)$, and no traveling waves exist for $c<c_*$. The uniqueness (up to shift) of the traveling waves was proved by
Aguerrea-Gomez-Trofimchuk \cite{Aguerrea-Gomez-Trofimchuk} by means of the Diekmann-Kaper theory. The wave $\phi(x+c_*t)$ with $c=c_*$ is called the critical traveling wave, and the wave with $c>c_*$ is usually called the non-critical traveling wave.  When $\frac{p}{\delta}>e$, the birth rate function $b(v)$ is non-monotone under consideration of $v\in [0,v_+]$, which causes that the equation \eqref{1.1} loses its monotonicity and the comparison principle does not hold. So the upper-lower solutions method cannot be applied to this case. By using the Lyapunov-Schmidt reduction method, Faria-Huang-Wu \cite{Faria-Huang-Wu} showed that, when the time-delay is small enough $r\ll 1$ and the wave speed is large enough $c\gg c_*$, the monotone  traveling waves $\phi(x+ct)$ exist. Later then, Faria-Trofimchuk \cite{Faria-Trofimchuk1} by analyzing heteroclinic solutions, and Ma \cite{Ma} by constructing auxiliary functions, both showed that, when $e<\frac{p}{\delta}\le e^2$, the traveling waves $\phi(x+ct)$ exist for all $c\ge c_*$. No restriction is on the time-delay $r>0$ in this case. But, when the time-delay $r$ is big: $r>\underline{r}>0$, where $\underline{r}>0$ given by
\begin{equation}\label{new-2017-3}
\delta\Big(\ln \frac{p}{\delta}-1\Big) \underline{r} e^{\delta \underline{r}+1}=1,
\end{equation}
then the traveling waves may occur oscillations. Here $\underline{r}$
 is the critical number for the time-delay $r$ such that the linear delay equation
 \begin{equation}\label{new-2017}
 v'(t)+\delta v(t)=b'(v_+)v(t-r)
 \end{equation}
  may have oscillating solutions.
Furthermore, when $\frac{p}{\delta}> e^2$, Gomez-Trofimchuk \cite{GT} proved that the  traveling waves exist only for $r<\bar{r}$
and no traveling waves for $r\ge \bar{r}$, where
\begin{equation}\label{new-2017-2}
\bar{r}:
=\frac{\pi-\arctan{\sqrt{(\ln\frac{p}{\delta}-2)\ln\frac{p}{\delta}}}}{\delta \sqrt{(\ln\frac{p}{\delta}-2)\ln\frac{p}{\delta}}}>0
\end{equation}
is the Hopf-bifurcation point to the equation \eqref{new-2017}.

In order to determine the minimal wave speed  $c_*>0$, let us linearize the equation \eqref{1.3} around $v_-=0$ for $\xi \sim -\infty$,
and test the eigenfunction by $\phi(\xi)=e^{\lambda\xi}$, we then have the characteristic equation for $\lambda>0$:
\[
c\lambda -D\lambda^2+\delta=pe^{-\lambda cr}.
\]
The left-hand-side functions $F_c(\lambda):=c\lambda -D\lambda^2+\delta$ and the right-hand-side function $G_c(\lambda):=pe^{-\lambda cr}$
of the above characteristic equation have a unique tangent point $(c_*,\lambda_*)$ with $c_*>0$ and $\lambda_*>0$, which is our minimal wave speed. For details, we refer to the graphes shown in \cite{Mei-Lin-Lin-So1}. Namely, $(c_*,\lambda_*)$
 is uniquely determined by
\begin{equation}
c_*\lambda_* -D\lambda_*^2  + \d=pe^{-\lambda_* c_*r} \ \mbox{ and } \
c_* -2D\lambda_*  =-c_*rpe^{-\lambda_* c_*r},
\label{1.6}
\end{equation}
and when $c>c_*$, there exist two numbers $\lambda_2>\lambda_1>0$ such that
\begin{equation}
c\lambda_i-D\lambda^2_i + \d=pe^{-\lambda_i cr}, \qquad \mbox{ for } i=1,2,
\label{nnew}
\end{equation}
and
\begin{equation}
c\lambda-D\lambda^2 + \d>pe^{-\lambda cr}, \qquad \mbox{ for } \lambda\in (\lambda_1,\lambda_2).
\label{1.6-2}
\end{equation}
As showed in \cite{Fang-Zhao,Faria-Trofimchuk1,GT,Ma,TT,TTT}, see also the summary in \cite{CMYZ,LLLM}, we have the following existence and uniqueness of the traveling waves, as well as the property of oscillations:

\begin{enumerate}
\item
When $e<\frac{p}{\d}\le e^2$, the traveling wave
$\phi(x+ct)$  exists  uniquely (up to a shift) for every $c\ge c_*=c_*(r)$, where the time-delay $r$ is  allowed to be any number in $[0,\infty)$.  If $0\le r < \underline{r}$, where $\underline{r}$, given by \eqref{new-2017-3},
is the critical point for the traveling waves to possibly  occur oscillations,
 then these traveling waves are monotone \cite{GT};   while, if the time delay $r \ge \underline{r}$, then the traveling waves are still monotone for
$(c,r)\in [c_*,c^*]\times [\underline{r},r_0]$, where $c_*=c_*(r)$ is the minimum wave speed as mentioned before, $c^*=c^*(r)$ is given by the characteristic equation for \eqref{1.3} around $v_+$, namely, the pair of $(c^*,\lambda^*)$ is determined by
\begin{equation}\label{c-up-star}
-c^*\lambda^*-D(\lambda^*)^2+d =b'(v_+)e^{\lambda^*c_*r},
\end{equation}
and $r_0$($>\underline{r}$) is the unique intersection point of two curves $c_*(r)$ and $c^*(r)$; and the traveling waves are oscillating around $v_+$ for $(c,r)\not\in [c_*,c^*]\times [\underline{r},r_0]$, namely, either $c>c^*$ or $r>r_0$ (c.f. \cite{GT,LLLM}).
\item
When $\frac{p}{\d}>e^2$,
on the other hand, the traveling wave $\phi(x+ct)$ with $c\ge  c_*$ can exist only  when $r< {\overline r}$,
and no traveling wave can exist  for  $r\ge {\overline r}$,
where ${\overline r}$ is the Hopf-bifurcation point given in \eqref{new-2017-2}, and the waves are monotone for $0 < r < \underline{r}$ and  oscillating for $ r \in (\underline{r},\overline{r})$ (c.f. \cite{GT,LLLM}).
\end{enumerate}

To investigate  the global stability of these oscillatory traveling waves will be our main target in this paper, including the critical wavefronts. The study on the critical traveling waves in the biological invasions is particularly interesting but also quite challenging,
 because the critical wave speed is usually the spreading speed for all solutions with
initial data having compact supports \cite{Liang-Zhao-1,Th-Zh}. In what follows,  we are going to review the progress on the stability of traveling waves for the type of mono-stable equations like Fisher-KPP equations without/with time-delay.

For regular mono-stable reaction-diffusion equations without time-delay ($r=0$) such as the classic Fisher-KPP equation, the existence of traveling waves and their stability have been one of the hot research spots. In 1976, by using the spectral analysis method,  Sattinger \cite{Sattinger} first proved
that, for given non-critical waves with $c>c_*$, when the initial perturbations around the waves are space-exponentially decay at the far field $-\infty$, then these non-critical waves are time-exponentially stable.
Since then, the study on stability of non-critical traveling waves has been intensively studied, for example, see \cite{AW1, Chen, Fife-McLeod, Ga, G-S, Hamel1,Hamel2,Hou-Li,Ki,Lau,M-R,S-S} and the references therein, see also the textbook \cite{VVV} and the survey paper \cite{Xin}.
However, the study on stability of critical traveling waves with $c=c_*$ is very limited, because this is critical case with  special difficulty. In 1978, by using the maximum principle method, Uchiyama \cite{Uchiyama} proved  the local stability for the traveling waves including
the critical waves, but no convergence rate
for the critical waves case was related. Later then,  Bramson \cite{Bramson} derived the sufficient and
necessary condition for the stability of noncritical and critical waves (no convergence
rates issued) by probability method.  Lau \cite{Lau} obtained the same results in a different way. Regarding the convergence rates to the critical traveling waves,
Moet \cite{Moet} first obtained the algebraic convergence rate $O(t^{-1/2})$ by using the Green function method.
 Kirchg\"assner \cite{Ki} then showed  the algebraic stability for the critical waves
in the form $O(t^{-1/4})$  by the spectral method. Gallay \cite{Ga} further improved the algebraic rate to
$O(t^{-3/2})$  by using the renormalization group method, when the
corresponding initial data converges   to the critical wave much fast like  $O(e^{-x^2/4})$ as $x\to-\infty$. More general case for parabolic equations was investigated by Eckmann and Wayne in \cite{E-W}.

For the mono-stable reaction-diffusion equations with time-delay, in 1987 Schaaf \cite{Schaaf} first proved the linear stability for the non-critical traveling waves by the spectral analysis method. This topic was not touched until in 2004 Mei-So-Li-Shen \cite{MSLS} proved the nonlinear stability of fast traveling waves with $c\gg c_*$ by the weighted energy method. When the equation is monotone (namely, $b(v)$ is increasing), Mei and his collaborators \cite{Mei-Lin-Lin-So1, Mei-Lin-Lin-So2, Mei-Ou-Zhao, Mei-So, Mei-Wang} further showed that all non-critical traveling waves are exponentially stable and all critical traveling waves are algebraically stable. When  $b(v)$ is non-monotone, the equation \eqref{1.1} losses its monotonicity. The solution is usually oscillating as the  time-delay $r$ is large, and  the traveling waves may occur oscillations around $v_+$, as theoretically proved in \cite{GT,TTT} and numerically reported in \cite{CMYZ,LLLM}. Such  oscillations   of the traveling waves are  interesting  and important  from both physical and mathematical points of view. Recently, by using the weighted energy method with the help of nonlinear Hanalay's inequality, Lin-Lin-Lin-Mei \cite{LLLM} proved that when the initial perturbation is small, then all non-critical oscillatory traveling waves are locally  stable with exponential convergence rate. Furthermore, by analyzing the decay rate of the (oscillatory) critical traveling waves at the unstable node $v_-=0$, and applying the anti-weighted energy method,  Chern-Mei-Yang-Zhang \cite{CMYZ} obtained  the local stability for the critical oscillatory traveling waves. But the convergence rate to the critical waves is still open.
 The interesting but also challenging questions are whether these oscillatory wavefronts are globally stable, and what will be  the optimal convergence rates, particularly the convergence rate for the critical wavefronts.  Note that,  the existing methods cannot be applied to our case, due to the lack of monotonicity of the equation and the waves, and the bad effect of time-delay. Of course, the critical wave case is always challenging as we know.

In this paper, the main targets are to prove that all critical or non-critical oscillatory traveling waves are globally stable,  and further to derive the optimal convergence rates to the critical/non-critical oscillatory traveling wave. That is, for all oscillatory traveling waves, including the critical traveling waves, the original solution of \eqref{1.1} will time-asymptotically converge to the targeted wave, even if the corresponding initial perturbation in a certain weighted space is big. The optimal convergence rates for the critical/non-critical  wavefronts are $O(t^{-\frac{1}{2}})$ and $O(t^{-\frac{1}{2}}e^{-\mu t})$, respectively.    As mentioned in \cite{CMYZ}, the usual approaches for deriving the  convergence rate are either the monotonic method with the help of  the decay estimates for linearized equations \cite{Mei-Lin-Lin-So1,Mei-Ou-Zhao}, or Fourier transform  \cite{Huang-Mei-Wang,Mei-Wang},  or the multiplier method \cite{RTY,Todorova}, or the method of approximate Green function \cite{Nishihara-Wang-Yang,Wang-Yang}. Since our equation is lack of monotonicity and the traveling waves are oscillating,  and the bad effect of time-delay, it seems that we could not be able to adopt these methods mentioned before.  However, we have some key observation, that is, although the equation and the waves are non-monotone, and we lose the comparison principle, due to the structure of the governing equation, we realize that, after using the anti-weight technique, the perturbed equation is reduced to a new equation, and  the absolute value of the oscillating solution for the new equation can be bounded by the positive solution of linear delayed  heat equation with constant coefficients. This can make us possibly to reach our goal by deriving the optimal decay estimate for the linearized solution, where, by  using Fourier transform, we derive the fundamental solution for the time-delayed linear heat equation and  its  optimal convergence rate. Based on the boundedness estimates for the oscillatory solution, we further prove the global stability for the oscillatory traveling waves. To our best knowledge, this is the first frame work to show the global stability for the oscillatory traveling waves with the optimal convergence rates, particularly, for the critical wave case.

Before stating our main results, we first introduce some notations on the solution spaces. Throughout this paper,  $C>0$ denotes for a generic
constant, and $C_i>0$ ($i=0,1,2,\cdots$) for specific positive constants.
$L^2(\R)$ is  the space of the square integrable functions, $H^k(\R)$ is the Sobolev space, $C(\R)$ is the space of bounded continuous functions,
$C([0,T] ; {\cal B})$ is the space of
the $\cal B$-valued continuous functions on $[0,T]$, where ${\cal B}$ is
a Banach space, and
 $T > 0$ is a number. Similarly,
$L^2([0,T] ; {\cal B})$ is the space of the ${\cal B}$-valued
$L^2$-functions on $[0,T]$.

To handle delay equation with delay $r$, as denoted in \cite{CMYZ,LLLM}, we define the uniformly continuous space   $\mathcal{C}_{unif}[-r,T]$, for $0<T\le \infty$,  by
\begin{align}
\mathcal{C}_{unif}[-r,T]:=
\{ &v(t,x) \in C([-r,T]\times \R) \mbox{ such that } \notag \\
 &\lim_{x\to+\infty} v(t,x) \mbox{ exists uniformly in } t\in[-r,T], \mbox{ and } \label{new4} \\
&
\lim_{x\to +\infty}v_x(t,x)=\lim_{x\to +\infty} v_{xx}(t,x)=0 \notag \\
& \mbox{ uniformly with respect to } t \in [-r,T] \}. \notag
\end{align}
For perturbed equation around the traveling waves, we now define the solution space. First, we introduce a weight function.
For $c\ge c_*$, we define a weight function
 \begin{equation}
 w(\xi) :=
 \begin{cases}
 e^{-2\lambda \xi}, \ \ \xi\in \R , & \mbox{ for } c>c_*, \ \lambda\in (\lambda_1,\lambda_2), \\
 e^{-2\lambda_* \xi}, \ \ \xi\in \R, & \mbox{ for } c=c_*,
 \end{cases}
 \label{new0-1}
 \end{equation}
 where $\xi=x+ct$ for $(t,x)\in[-r,\infty)\times \R$, and the number $\lambda$ for the case of $c>c_*$ is selected between  $\lambda_1$ and $\lambda_2$, and
 $\lambda_1$ and $\lambda_2$ are specified in \eqref{nnew}.
Notice that, for $c\ge c_*$,  $\lim_{\xi\to-\infty}w(\xi)=+\infty$ and $\lim_{\xi\to+\infty}w(\xi)=0$, because $\lambda>0$ and $\lambda_* > 0$.
We also denote the weighted Sobolev space $W^{2,1}_w(\R)$ by
\[
W^{2,1}_w(\R)=\{u|wu\in L^1(\R), w\partial^i_\xi u \in L^1(\R), \ i=1,2\}.
\]

 Our main stability theorems are as follows.

 \begin{theorem}[Global stability with optimal convergence rates]\label{rates}
 For any given  traveling wave $\phi(x+ct)=\phi(\xi)$ with $c\ge c_*$, no matter it is oscillatory or not, when the initial perturbation satisfies  $v_0(s,\xi)-\phi(\xi)\in C_{unif}(-r,0)\cap C([-r,0];W^{2,1}_w(\R))$ and $\partial_s(v_0-\phi)\in L^1([-r,0];L^1_w(\R))$,  then the following global stability holds.
 \begin{enumerate}
 \item When  $e<\frac{p}{\d}\le e^2$,  for any time-delay $r>0$,
then
\begin{equation}
\begin{cases}
\sup_{x\in \R} |v(t,x)-\phi(x+ct)|\le C t^{-\frac{1}{2}}e^{-\mu t} \ \ &\mbox{ for } c>c_*, \\
\sup_{x\in \R} |v(t,x)-\phi(x+c_*t)|\le C t^{-\frac{1}{2}} \ \ &\mbox{ for } c=c_*,
\end{cases}
\label{2345}
\end{equation}
where $\mu$ is  a positive number  satisfying
\begin{equation}\label{1234}
0<\mu<\min\{ \d, c\lambda+\d-D\lambda^2-pe^{-\lambda cr}\}, \ \ \lambda\in(\lambda_1,\lambda_2).
\end{equation}

 \item When $\frac{p}{\d}>e^2$  but with a small time-delay $0<r<{\overline r}$, where $\bar{r}$ is defined in \eqref{new-2017-2},
 then the stability \eqref{2345} holds.

\end{enumerate}
\end{theorem}

\medskip

\begin{remark}
\begin{enumerate}
\item When $1<\frac{p}{\delta}\le e$, the traveling waves $\phi(x+ct)$ for all $c\ge c_*$ are monotone, and the global stability of these monotone critical/non-critical waves as well as the optimal convergence rates were shown by Mei-Ou-Zhao in \cite{Mei-Ou-Zhao}.
\item When $e<\frac{p}{\delta}\le e^2$ with any size of time-delay $r>0$, or $\frac{p}{\delta}>e^2$ but with small time-delay $r\in (0,\bar{r})$,   the local stability for the critical/non-critical wavefronts  was proved in \cite{Chern-Mei-Zhang-Zhang, Lin-Lin-Lin-Mei}, respectively. Here in Theorem \ref{rates}, we obtain  the global stability of these oscillatory traveling waves, where the initial perturbations around the traveling waves in the weighted space $C_{unif}(-r,0)\cap C([-r,0];W^{2,1}_w(\R))$ can be arbitrarily large. To our best knowledge, this is the first work to show the global stability for the oscillatory waves.

\item When $c=c^*$, we obtain the optimally algebraic convergence rate $O(t^{-\frac{1}{2}})$ to the critical oscillating waves. This answers the open question left in \cite{CMYZ}. The algebraic convergence rate also matches what shown for the critical monotone traveling waves for Fisher-KPP equations with or without time-delay \cite{Mei-Ou-Zhao,Mei-Wang,Moet}.

\item When $c>c^*$, the obtained exponential convergence $O(e^{-\mu}t^{-\frac{1}{2}})$ improves  the previous study for the non-critical oscillating waves in \cite{LLLM}.

\end{enumerate}
\end{remark}

\medskip

\section{Proof of main theorem}

Let $\phi(x+ct)=\phi(\xi)$ be any given monotne/non-monotone traveling wave with $c\ge c_*$, and define
\[
u(t,\xi):=v(t,x)-\phi(x+ct), \ \ u_0(s,\xi):=v_0(s,x)-\phi(x+cs), \ \ c\ge c_*.
\]
Then, from \eqref{1.1}, $u(t,\xi)$ satisfies, for $c\ge c_*$,
\begin{equation}
\begin{cases}
 \dfrac{\partial u}{\partial t} +c\dfrac{\partial u}{\partial\xi}- D\dfrac{\partial^2
u}{\partial \xi^2} +\d u    = P(u(t-r,\xi-cr)),  \ \ (t,\xi)\in
\R_+\times \R,\\
u(s,\xi)=u_0(s,\xi), \ \ s\in [-r,0], \ \xi\in \R,
\end{cases}
\label{2.1}
\end{equation}
where
\begin{equation}
P(u):=b(\phi+u)-b(\phi)=b'(\tilde\phi)u,
\label{2.2}
\end{equation}
for some $\tilde\phi$ between $\phi$ and $\phi+u$, with $\phi=\phi(\xi-cr)$ and $u=u(t-r,\xi-cr)$ for $c\ge c_*$.

We first prove the existence and uniqueness of solution to the initial value problem \eqref{2.1} in the uniformly continuous space $C_{unif}[-r,\infty)$.

\begin{proposition}[Existence and Uniqueness] When $e<\frac{p}{\d}\le e^2$ for any $r>0$, or $\frac{p}{\d}>e^2$ but for $r\in (0,\bar{r})$,  if the
 initial perturbation $u_0\in C_{unif}[-r,0]$ for $c\ge c_*$,
 then
 the solution $u(t,\xi)$ of the perturbed equation \eqref{2.1}  is unique and time-globally exits in $C_{unif}[-r,\infty)$.
\end{proposition}
{\bf Proof.} When $t\in[0,r]$, since $t-r\in [-r,0]$ and $u(t-r,\xi-cr)=u_0(t-r,\xi-cr)$, then \eqref{2.1} reduce to the following linear equation
\begin{equation}
\begin{cases}
 \dfrac{\partial u}{\partial t} +c\dfrac{\partial u}{\partial\xi}- D\dfrac{\partial^2
u}{\partial \xi^2} +\d u    = P(u_0(t-r,\xi-cr)),  \ \ (t,\xi)\in
\R_+\times \R,\\
u(0,\xi)=u_0(0,\xi), \ \ \ \xi\in \R,
\end{cases}
\label{2.1-new-1}
\end{equation}
and it possesses a unique solution $u(t,\xi)$ for $t\in [0,r]$ in the integral form of
\begin{eqnarray}
u(t,\xi)&=&e^{-\d t} \int_{-\infty}^\infty G(t,\eta)u_0(0,\xi-\eta)d\eta \notag \\
&&+\int^t_0 e^{-\d (t-s)} \int_{-\infty}^\infty G(t-s,\eta)P(u_0(s-r,\xi-\eta-cr)d\eta ds,
\label{2.1-new-2}
\end{eqnarray}
where $G(t,\eta)$ is the heat kernel
\[
G(t,\eta)=\frac{1}{\sqrt{4\pi Dt}}e^{-\frac{(\eta+ct)^2}{4Dt}}, \ \ c\ge c_*.
\]
Since $u_0\in C_{unif}[-r,0]$, namely, $\lim_{\xi\to\infty}  u_0(t,\xi)=u_0(t,\infty)$ and $\lim_{\xi\to\infty}  u_{0,\xi}(t,\xi)=0$ and
$\lim_{\xi\to\infty}  u_{0,\xi\xi}(t,\xi)=0$ uniformly in $t\in [-r,0]$, we immediately prove the following uniform convergence
\begin{eqnarray}
\lim_{\xi\to\infty} u(t,\xi)&=&e^{-\d t} \int_{-\infty}^\infty G(t,\eta)\lim_{\xi\to\infty}u_0(0,\xi-\eta)d\eta \notag \\
&&+\int^t_0 e^{-\d (t-s)} \int_{-\infty}^\infty G(t-s,\eta)\lim_{\xi\to\infty}P(u_0(s-r,\xi-\eta-cr)d\eta ds\notag \\
&=&e^{-\d t} u_0(0,\infty)\int_{-\infty}^\infty G(t,\eta)d\eta \notag \\
&&+\int^t_0 e^{-\d(t-s)} P(u_0(s-r,\infty) \int_{-\infty}^\infty G(t-s,\eta)d\eta ds \notag \\
&=&e^{-\d t} u_0(0,\infty)
+\int^t_0 e^{-\d(t-s)} P(u_0(s-r,\infty) ds\notag \\
&=:&g_1(t), \ \ \mbox{ uniformly in } t\in [0,r],
\label{2.1-new-3}
\end{eqnarray}
and
\begin{eqnarray}
\lim_{\xi\to\infty} \partial_\xi^k u(t,\xi)&=&e^{-\d t} \int_{-\infty}^\infty \partial_\eta^k G(t,\eta)\lim_{\xi\to\infty}u_0(0,\xi-\eta)d\eta \notag \\
&&+\int^t_0 e^{-\d(t-s)} \int_{-\infty}^\infty \partial_\eta^k G(t-s,\eta)\lim_{\xi\to\infty}P(u_0(s-r,\xi-\eta-cr)d\eta ds\notag \\
&=&e^{-\d t} u_0(0,\infty)\int_{-\infty}^\infty \partial_\eta^k G(t,\eta)d\eta \notag \\
&&+\int^t_0 e^{-\d(t-s)} P(u_0(s-r,\infty) \int_{-\infty}^\infty  \partial_\eta^k G(t-s,\eta)d\eta ds \notag \\
&=&0, \ \ \mbox{ for } k=1,2, \  \mbox{ uniformly in } t\in [0,r].
\label{2.1-new-4}
\end{eqnarray}
Here, we used the facts
\[
\int_{-\infty}^\infty  G(t-s,\eta)d\eta=1, \ \ \   G(t-s,\pm\infty)=0, \ \ \mbox{ and } \partial_\eta G(t-s,\eta)\Big|_{\eta=\pm\infty}=0.
\]
Thus, we have proved $u\in C_{unif}[-r,r]$.

Now we consider \eqref{2.1} for $t\in [r,2r]$. Since $t-r \in [0,r]$ and $u(t-r,\xi-cr)$ is solved already in last step for \eqref{2.1-new-1}, thus $P(u(t-r,\xi-cr))$ is known for \eqref{2.1} with  $t\in [0,2r]$, namely,  the equation \eqref{2.1} is linear for $t\in [0,2r]$. As showed before, we can similarly prove  the existence and uniqueness of the solution $u$ to \eqref{2.1} for $t\in [0,2r]$, and particularly
 $u\in C_{unif}[-r,2r]$.

By repeating this procedure for $t\in [nr,(n+1)r]$ with $n\in {\mathbb{Z}}_+$ (the set of all positive integers), we prove that there exists a unique solution $u\in C_{unif}([-r,(n+1)r]$ for \eqref{2.1}, and step by step, we finally prove the uniqueness and time-global existence of the solution $u\in C_{unif}[-r,\infty)$ for \eqref{2.1}. The proof is complete.
$\square$

\medskip

The most important part of the paper is to prove the following global stability with the optimal convergence rates.

\begin{proposition}[Stability with optimal convergence rates]\label{stability-2} When $e<\frac{p}{\d}\le e^2$ for any $r>0$, or $\frac{p}{\d}>e^2$ but for $r\in (0,\bar{r})$, if  $ u_0\in C_{unif}(-r,0)\cap  L^1([-r,0];W^{2,1}_w(\R))$ and $\partial_s u_0 \in L^1([-r,0];L^{1}_w(\R))$ for $c\ge c_*$, no matter how large the initial perturbation $u_0$ is, then
\begin{itemize}
\item
when $c=c_*$, it holds
\begin{equation}
\sup_{\xi\in \R} |u(t,\xi)|\le C t^{-\frac{1}{2}}.
\label{alg-stable-2}
\end{equation}
\item when $c>c_*$, it holds
\begin{equation}
\sup_{\xi\in \R} |u(t,\xi)|\le C t^{-\frac{1}{2}}e^{-\mu t}
\label{exp-stable-2}
\end{equation}
for a positive constant $\mu$ specified in \eqref{1234}.
\end{itemize}

\end{proposition}

In order to prove Proposition \ref{stability-2}, we need serval lemmas to complete it. Since $b'(0)=p>\d=d'(0)$, namely, $u_-=0$ is the unstable node of
\eqref{2.1}, heuristically, for a general initial data $u_0$, we cannot expect the convergence $u\to 0$ as $t\to\infty$. But, inspired by
\cite{Schaaf,Sattinger} for the equation \eqref{2.1} with or without time-delay, we expect the solution $u$ to decay to zero when the initial perturbation is only exponentially decay at the far field $\xi=-\infty$. We need the anti-weighted energy technique to treat this problem.
Let us define
\begin{equation}
u(t,\xi)=[w(\xi)]^{-\frac{1}{2}}\tilde{u}(t,\xi), \ \ i.e., \ \ \tilde{u}(t,\xi)=\sqrt{w(\xi)}u(t,\xi)=e^{-\lambda \xi}u(t,\xi),
\end{equation}
where $\lambda \in (\lambda_1,\lambda_2)$ for $c>c_*$ and $\lambda=\lambda_*$ for $c=c_*$,
we get the following equations for the new unknown $\tilde{u}(t,\xi)$ with $c\ge c_*$:
\begin{equation}
\begin{cases}
 \dfrac{\partial \tilde{u}}{\partial t} - D\dfrac{\partial^2
\tilde{u}}{\partial \xi^2} + a_0(c) \dfrac{\partial\tilde{u} }{\partial \xi} + a_1(c) \tilde{u}
 = \tilde{P}(\tilde{u}(t-r,\xi-cr)),  \ \ (t,\xi)\in
\R_+\times \R,\\
\tilde{u}(s,\xi)=\sqrt{w(\xi)}u(s,\xi)=:\tilde{u}_0(s,\xi), \ \ s\in [-r,0], \ \xi\in \R,
\end{cases} \label{2.2-new}
\end{equation}
where
\begin{equation}
a_0(c):=c-2D\lambda, \ \ \ a_1(c):=c\lambda+\d-D\lambda^2, \ \mbox{ for } c\ge c_*,
\label{2.3-new}
\end{equation}
which satisfy (by \eqref{1.6} and \eqref{1.6-2})
\begin{eqnarray}
&a_1(c)=c\lambda +\d-D\lambda^2>pe^{-\lambda cr}, &\mbox{ for } c>c_*, \lambda\in(\lambda_1,\lambda_2), \label{04-28-1}\\
&a_1(c_*)=c_*\lambda_*+\d-D\lambda_*^2=pe^{-\lambda_* c_*r}, &\mbox{ for } c=c_*.  \label{04-28-2}
\end{eqnarray}
Here,
\begin{equation}
\tilde{P}(\tilde{u})=e^{-\lambda \xi}P(u)
\label{2.4-new}
\end{equation}
satisfies (by Taylor's expansion formula)
\begin{eqnarray}
\tilde{P}(\tilde{u}(t-r,\xi-cr))&=& e^{-\lambda \xi}P(u(t-r,\xi-cr))\notag \\
&=& e^{-\lambda \xi}b'(\tilde{\phi})u(t-r,\xi-cr)\notag \\
&=&e^{-\lambda c r}b'(\tilde{\phi}) \tilde{u}(t-r,\xi-cr), \ \ c\ge c_*,
\label{2.4-new-1}
\end{eqnarray}
for some function $\tilde{\phi}$ between $\phi$ and $\phi+u$ with $c\ge c_*$ (see \eqref{2.2}).
We further have
\begin{equation}
|\tilde{P}(\tilde{u}(t-r,\xi-cr))|\le  p e^{-\lambda cr}|\tilde{u}(t-r,\xi-cr)|.
\label{2.4-new-1-new}
\end{equation}

Since $b'(s)$ can be negative for $s\in (0,v_+)$, then the solution $\tilde u$ for the equation \eqref{2.2-new} with the nonlinear term \eqref{2.4-new-1} will be  oscillating around $v_+$ when the time-delay $r$ is large, as numerically reported in \cite{CMYZ,LLLM}, and the comparison principle doesn't hold in this case. So the monotonic technique cannot be applied. On the other hand, in \eqref{2.4-new-1}, since the coefficient $b'(\tilde{\phi})$ is variable, we are unable to derive the decay rate directly by applying Fourier's transform. However, by a deep observation we may  establish a crucial boundedness estimate for the oscillating solution ${\tilde u}(t,\xi)$. In order to look for such a boundedness, let us consider the following linear delayed reaction-diffusion equation, for $c\ge c_*$,
\begin{equation}
\begin{cases}
 \dfrac{\partial u^+}{\partial t} - D\dfrac{\partial^2
u^+}{\partial \xi^2} + a_0(c) \dfrac{\partial u^+ }{\partial \xi} + a_1(c) u^+ \\
  \qquad = pe^{-\lambda cr}u^+(t-r,\xi-cr)),  \ \ (t,\xi)\in
\R_+\times \R,\\
u^+(s,\xi)=u^+_0(s,\xi)\ge 0, \ \ s\in [-r,0], \ \xi\in \R.
\end{cases} \label{2.14}
\end{equation}

\begin{lemma} When $u_0^+(s,\xi)\ge 0$ for $(s,\xi)\in [-r,0]\times \R$, then
$u^+(t,\xi)\ge 0$ for $(t,\xi)\in [-r,\infty)\times \R$.
\end{lemma}
{\bf Proof}. As showed in \cite{Mei-Lin-Lin-So1}, we can similarly prove the positiveness of the solution $u^+(t,\xi)$.
In fact, when $t\in [0,r]$, we have $t-r\in [-r,0]$ and
\[pe^{-\lambda cr}u^+(t-r,\xi-cr))=pe^{-\lambda cr}u^+_0(t-r,\xi-cr))\ge 0.
\]
Then \eqref{2.14} implies
\[
\dfrac{\partial u^+}{\partial t} - D\dfrac{\partial^2
u^+}{\partial \xi^2} + a_0(c) \dfrac{\partial u^+ }{\partial \xi} + a_1(c) u^+\ge 0,
\]
which gives $u^+(t,\xi)\ge 0$ for $t\in[0,r]$. Repeating this procedure step by step, we can prove $u^+(t,\xi)\ge 0$ for $t\in[nr,(n+1)r]$, and further
$u^+(t,\xi)\ge 0$ for $t\in R_+$. The proof is complete. $\square$

Next we establish the following crucial boundedness estimate for the oscillating solution ${\tilde u}(t,\xi)$ of \eqref{2.1}.

\begin{lemma}[Boundedness estimate]\label{lem2.1} Let ${\tilde u}(t,\xi)$ and $u^+(t,\xi)$ be the solutions of \eqref{2.1} and \eqref{2.14}, respectively. When
\begin{equation}
|{\tilde u}_0 (s,\xi)|\le u^+_0(s,\xi), \ \ \ (s,\xi)\in [-r,0]\times \R,
\label{2.16}
\end{equation}
then
\begin{equation}
|{\tilde u} (t,\xi)|\le u^+(t,\xi), \ \ \ (t,\xi)\in \R_+\times \R.
\label{2.17}
\end{equation}
\end{lemma}
{\bf Proof.} First of all, we prove $|{\tilde u}(t,\xi)|\le u^+(t,x)$ for $t\in [0,r]$. In fact, when $t\in [0,r]$, namely, $t-r\in [-r,0]$,
from \eqref{2.16}, we have
\begin{eqnarray}
|{\tilde u}(t-r,\xi-cr)|&=&|{\tilde u}_0(t-r,\xi-cr)|\notag \\
&\le& u^+_0(t-r,\xi-cr)\notag \\
&=&u^+(t-r,\xi-cr), \ \ \mbox{ for } t\in[0,r] \mbox{ and } c\ge c_*.
\label{2.18}
\end{eqnarray}
Let
\[
U^-(t,\xi):=u^+(t,\xi) - {\tilde u}(t,\xi) \ \mbox{ and }  \ U^+(t,\xi):=u^+(t,\xi)+{\tilde u}(t,\xi),
\]
we are going to estimate $U^\pm(t,\xi)$ respectively.

From \eqref{2.2-new}, \eqref{2.4-new-1} and \eqref{2.14}, then $U^-(t,\xi)$ satisfies, for $c\ge c_*$,
\begin{eqnarray}
& & \dfrac{\partial U^-}{\partial t} +a_0(c)\dfrac{\partial U^-}{\partial\xi}- D\dfrac{\partial^2
U^-}{\partial \xi^2} + a_1(c)  U^- \notag \\
 & &= b'(0)e^{-\lambda c r} u^+(t-r,\xi-cr)-b'(\tilde{\phi}) e^{-\lambda c r} {\tilde u}(t-r,\xi-cr) \notag \\
 & &\ge  b'(0)e^{-\lambda c r} u^+(t-r,\xi-cr)-| b'(\tilde{\phi})| e^{-\lambda c r} |{\tilde u}(t-r,\xi-cr)| \notag \\
 & &\ge 0, \ \ t\in [0,r],
\label{2.19}
\end{eqnarray}
where $\lambda \in (\lambda_1,\lambda_2)$ for $c>c_*$ and $\lambda=\lambda_*$ for $c=c_*$. Here we used 
\eqref{2.18} for $|\tilde u(t-r,\xi-cr)|\le |u^+(t-r,\xi-cr)|$ and the fact
\begin{equation}\label{new-2017-4}
|b'(\tilde{\phi})|=p|1-a\tilde{\phi}|e^{-a\tilde{\phi}} \le p=b'(0) \ \ \mbox{ for any } \tilde{\phi}>0.
 \end{equation}
  Thus,
\eqref{2.19} with the initial data $U^-_0(s,\xi)=u^+_0(s,\xi)-{\tilde u}_0(s,\xi)\ge 0$ reduces  to
\[
\begin{cases}
\dfrac{\partial U^-}{\partial t} +a_0(c)\dfrac{\partial U^-}{\partial\xi} - D\dfrac{\partial^2
U^-}{\partial \xi^2} + a_1(c) U^-\ge 0, \\
U^-_0(0,\xi)\ge 0,
\end{cases}
\]
which, by the regular comparison principle for parabolic equation without delay,  implies
\begin{equation}
U^-(t,\xi)=u^+(t,\xi)-{\tilde u} (t,\xi)\ge 0 \ \mbox{ for } (t,\xi)\in [0,r]\times \R \ \mbox{ and } c\ge c_*.
\label{2.20}
\end{equation}
On the other hand, $U^+(t,\xi)$ satisfies
\begin{eqnarray}
& & \dfrac{\partial U^+}{\partial t} +a_0(c)\dfrac{\partial U^+}{\partial\xi}- D\dfrac{\partial^2
U^+}{\partial \xi^2} + a_1(c) U^+ \notag \\
 & &= b'(0)e^{-\lambda c r} u^+(t-r,\xi-cr)+ b'(\tilde{\phi})e^{-\lambda c r}  {\tilde u}(t-r,\xi-cr) \notag \\
 & &\ge  b'(0)e^{-\lambda c r} u^+(t-r,\xi-cr)-| b'(\tilde{\phi})| e^{-\lambda c r} |{\tilde u}(t-r,\xi-cr)| \notag \\
 & &\ge 0, \ \ t\in [0,r],
\label{2.21}
\end{eqnarray}
because of $|b'(\tilde{\phi})|<b'(0)$ as mentioned in \eqref{new-2017-4},  and $|\tilde u(t-r,\xi-cr)|\le |u^+(t-r,\xi-cr)|$ from \eqref{2.18}. Therefore, we can similarly prove that
 \eqref{2.21} with the initial data $U^+_0(s,\xi)=u^+_0(s,\xi)+{\tilde u}_0(s,\xi)\ge 0$ implies
\begin{equation}
U^+(t,\xi)=u^+(t,\xi)+{\tilde u}(t,\xi)\ge 0 \ \mbox{ for } (t,\xi)\in [0,r]\times \R \ \mbox{ and } c\ge c_*.
\label{2.22}
\end{equation}
Combining \eqref{2.20} and \eqref{2.22}, we prove
\begin{equation}
|{\tilde u}(t,\xi)|\le u^+(t,\xi) \ \mbox{ for } (t,\xi)\in [0,r]\times \R \ \mbox{ and } c\ge c_*.
\label{2.23}
\end{equation}
Next, when $t\in [r,2r]$, namely, $t-r\in [0,r]$, based on \eqref{2.23} we can similarly prove
\[
\begin{cases}
U^-(t,\xi)=u^+(t,\xi)-{\tilde u}(t,\xi)\ge 0,\\
U^+(t,\xi)=u^+(t,\xi)+{\tilde u}(t,\xi)\ge 0,
\end{cases}
\ \mbox{ for } (t,\xi)\in [r,2r]\times \R \ \mbox{ and } c\ge c_*,
\]
namely,
\begin{equation}
|{\tilde u}(t,\xi)|\le u^+(t,\xi) \ \mbox{ for } (t,\xi)\in [r,2r]\times \R \ \mbox{ and } c\ge c_*.
\label{2.24}
\end{equation}
Repeating this procedure, we then further prove
\begin{equation}
|{\tilde u}(t,\xi)|\le u^+(t,\xi) \ \mbox{ for } (t,\xi)\in [nr,(n+1)r]\times \R \ \mbox{ and } c\ge c_*, \ \ n=1,2,\cdots,
\label{2.25}
\end{equation}
which implies
\begin{equation}
|{\tilde u}(t,\xi)|\le u^+(t,\xi) \ \mbox{ for } (t,\xi)\in [0,\infty)\times \R \ \mbox{ and } c\ge c_*.
\label{2.25-2}
\end{equation}
The proof is complete. $\square$

Next we derive the global stability with the optimal convergence rates for the linear equation \eqref{2.14}  by using the weighted energy method and
by carrying out the crucial estimates on the fundamental solutions. In order to derive the optimal decay rates for the solution of \eqref{2.14} in the cases of $c>c_*$ and $c=c_*$, respectively, we  first need to derive the fundamental solution, then show  the optimal decay rates of the fundamental solution.
Now let us recall the properties of the solutions to the delayed ODE.

\begin{lemma}[\cite{KIK}]
Let $z(t)$ be the solution to the following linear time-delayed ODE with time-delay $r>0$ and two constants $k_1$ and $k_2$
\begin{equation}\label{p1}
\begin{cases}
\displaystyle\frac{d}{dt}z(t)+k_1 z(t) =k_2 z(t-r) \\
z(s)=z_0(s), \ \ \ s\in[-\tau,0].
\end{cases}
\end{equation}
Then
\begin{equation}
z(t)=e^{-k_1(t+r)} e^{{\bar k_2}t}_\tau z_0(-r) +\int^0_{-r} e^{-k_1(t-s)}e^{{\bar k_2}(t-r-s)}_\tau [z_0'(s)+k_1 z_0(s)]ds ,
\label{p2}
\end{equation}
where
\begin{equation}
{\bar k_2}:=k_2 e^{k_1 \tau}, \label{k_2}
\end{equation}
and
$e^{{\bar k_2}t}_\tau$ is the so-called {\tt delayed exponential function} in the form
\begin{equation}
e^{{\bar k_2} t}_r =\begin{cases}
0, & -\infty<t<-r, \\
1, & -r\le t<0, \\
1+\frac{{\bar k_2} t}{1!}, & 0\le t<r, \\
1+\frac{{\bar k_2} t}{1!}+\frac{{\bar k_2}^2(t-r)^2}{2!}, & r\le t<2r, \\
\vdots & \vdots \\
1+\frac{{\bar k_2} t}{1!}+\frac{{\bar k_2}^2(t-r)^2}{2!}+\cdots + \frac{{\bar k_2}^m[t-(m-1)r]^m}{m!}, & (m-1)r\le t<mr, \\
\vdots & \vdots
\end{cases}
\label{p3}
\end{equation}
and $e^{{\bar k_2} t}_r$ is the fundamental solution to
\begin{equation}\label{p4}
\begin{cases}
\displaystyle\frac{d}{dt}z(t) ={\bar k_2} z(t-r) \\
z(s)\equiv 1, \ \ \ s\in[-r,0].
\end{cases}
\end{equation}
\end{lemma}

The property of the solution to the delayed linear ODE \eqref{p1} is well-known \cite{Mei-Wang}.

\begin{lemma}[\cite{Mei-Wang}]\label{lemma3} Let $k_1\ge 0$ and $k_2\ge 0$. Then the solution $z(t)$  to \eqref{p1} (or equivalently \eqref{p2}) satisfies
\begin{equation}
|z(t)|\le C_0 e^{-k_1 t} e^{{\bar k_2}t}_r, \label{12-26-1}
\end{equation}
where
\begin{equation}
C_0:= e^{-k_1\tau}|z_0(-r)| + \int^0_{-r} e^{k_1 s}|z'_0(s)+k_1 z_0(s)| ds, \label{12-26-2}
\end{equation}
and the fundamental solution $e^{{\bar k_2} t}_r$ with ${\bar k_2}>0$ to \eqref{p4} satisfies
\begin{equation}
e^{{\bar k_2} t}_\tau \le C(1+t)^{-\gamma} e^{{\bar k_2} t},   \ \ t>0
\label{p5}
\end{equation}
for arbitrary number $\gamma>0$.

Furthermore, when $k_1\ge k_2\ge 0$,  there exists a  constant $0< \varepsilon_1 <1$
such that
\begin{equation}
e^{-k_1 t}e^{{\bar k_2}t}_\tau \le C e^{-\varepsilon_1(k_1-k_2) t}, \ \ t>0
\label{12-26-3}
\end{equation}
and the solution $z(t)$ to \eqref{p1} satisfies
\begin{equation}
|z(t)|\le C e^{-\varepsilon_1(k_1-k_2) t}, \ \ t>0.
\label{12-26-4}
\end{equation}
\end{lemma}

Let us take Fourier transform to \eqref{2.14}, and denote the Fourier transform of $u^+(t,\xi)$ by $\hat{u}^+(t,\eta)$, that is,
\begin{equation} \label{07-29-1}
\begin{cases}
\displaystyle{\frac{d}{d t}\hat{u}^+(t,\eta) + A(\eta) \hat{u}^+(t,\eta) = B(\eta) \hat{u}^+(t-r,\eta),} \\
\hat{u}^+(s,\eta)=\hat{u}^+_0(s,\eta), \  s\in [-r,0], \ \ \eta\in R,
\end{cases}
\end{equation}
where
\begin{equation} \label{07-29-2}
\begin{cases}
A(\eta):= D|\eta|^2 + a_1(c) + \mbox{i} a_0(c)\eta, \\
B(\eta):=b'(0)e^{\mbox{i}cr\eta} e^{-\lambda cr},
\end{cases}
\mbox{ for } c\ge c_*.
\end{equation}
From  \eqref{p2},  the  linear time-delayed ordinary differential equation \eqref{07-29-1} can be solved by
\begin{eqnarray}
\hat{u}^+(t,\eta)&=& e^{-A(\eta)(t+r)} e_r^{\bar{B}(\eta)t} \hat{u}^+_0(-r,\eta) \notag \\
& & + \int^0_{-r} e^{-A(\eta)(t-s)}e_r^{\bar{B}(\eta)(t-r-s)} \notag \\
& & \qquad \ \ \times \Big[\frac{d}{ds}\hat{u}^+_0(s,\eta) + A(\eta)\hat{u}^+_0(s,\eta)\Big] ds,
\label{07-29-3}
\end{eqnarray}
where
\begin{equation}
\bar{B}(\eta):=B(\eta) e^{A(\eta)r},
\label{07-29-4}
\end{equation}

Taking the inverse Fourier transform to \eqref{07-29-3}, we have
\begin{eqnarray}
u^+(t,\xi)&=& \frac{1}{2\pi} \int^\infty_{-\infty}e^{\mbox{i}\xi\eta}e^{-A(\eta)(t+r)}e_r^{\bar{B}(\eta)t} \hat{u}^+_0(-r,\eta)d\eta \notag \\
& & + \frac{1}{2\pi} \int^0_{-r}  \int^\infty_{-\infty}e^{\mbox{i}\xi\eta}e^{-A(\eta)(t-s)}e_r^{\bar{B}(\eta)(t-r-s)} \notag \\
& & \qquad \qquad  \times \Big[\frac{d}{ds}\hat{u}^+_0(s,\eta) + A(\eta)\hat{u}^+_0(s,\eta)\Big] d\eta ds.
\label{07-29-5}
\end{eqnarray}

Next, we are going to estimate the decay rates for the solution $u^+(t,\xi)$.

\begin{lemma}[Optimal decay rates for linear delayed  equation]\label{lem2.2} Let the initial data $u_0^+(s,\xi)$ be  such that $u_0^+\in L^1([-r,0];W^{2,1}(\R))$ and $\partial_s u^+_0 \in L^1([-r,0];L^1(\R))$, then
\begin{equation}
\|u^+(t)\|_{L^\infty(\R)}\le
\begin{cases}
C(1+t)^{-\frac{1}{2}}e^{-\mu_1 t}, & \mbox{ for } c>c_*,\\
C(1+t)^{-\frac{1}{2}}, & \mbox{ for } c=c_*,
\end{cases}
\label{2.26}
\end{equation}
for some $0<\mu_1<c\lambda+\d-c\lambda^2-pe^{-c\lambda r}$ for $c>c_*$.

\end{lemma}
{\bf Proof}. Using Parseval's inequality, from \eqref{07-29-3} we have
\begin{eqnarray}
\|u^+(t)\|_{L^\infty(\R)}&\le&\|\hat{u}^+(t)\|_{L^1(\R)} \notag \\
&\le & \int_{-\infty}^\infty \Big|e^{-A(\eta)(t+r)}e_r^{\bar{B}(\eta)t} \hat{u}^+_0(-r,\eta)\Big|d\eta \notag \\
& & + \int^0_{-r} \int_{-\infty}^\infty \Big| e^{-A(\eta)(t-s)}e_r^{\bar{B}(\eta)(t-r-s)} \notag \\
& & \qquad \qquad   \times \Big[\frac{d}{ds}\hat{u}^+_0(s,\eta) + A(\eta)\hat{u}^+_0(s,\eta)\Big] \Big| d\eta ds \notag \\
&=:& I_1(t) +I_2(t).
\label{07-29-6}
\end{eqnarray}
To estimate $I_i(t)$ for $i=1,2$, from the definitions of $\bar{B}(\eta)$ (see \eqref{07-29-4} and \eqref{07-29-2}) and the delayed-exponential function $e^{{\bar k}_2t}_r$ (see \eqref{p3}), we first note
\begin{equation}\label{8-1-1}
|e^{-A(\eta)(t+r)}|= e^{-(D\eta^2+a_1(c))(t+r)}=e^{-k_1(c,\eta)(t+r)}, \ \ \mbox{ for } c\ge c_*,
\end{equation}
with
\begin{equation}\label{8-1-1-1}
k_1(c,\eta):=D\eta^2+a_1(c),
\end{equation}
and
\begin{eqnarray}\label{8-1-2}
|\bar{B}(\eta)|&\le& |B(\eta)|e^{|A(\eta)|r}
\le b'(0)e^{-\lambda cr} e^{k_1(c,\eta)} \notag \\
&=:& {\bar k}_2(c,\eta), \ \ \mbox{ for } c\ge c_*,
\end{eqnarray}
where
\begin{equation}\label{9-1-2-1}
{\bar k}_2(c,\eta):=k_2(c)e^{k_1(c,\eta)}, \ \ \mbox{and } k_2(c):= b'(0)e^{-\lambda cr},
\end{equation}
and
\begin{equation}\label{8-1-3}
|e^{{\bar B}(\eta)t}_r| \le e^{|{\bar B}(\eta)|t}_r = e^{{\bar k}_2(c,\eta)t}_r, \ \ \mbox{ for } c\ge c_*.
\end{equation}
Noting \eqref{1.6} and \eqref{1.6-2}, we have
\[
k_1(c,\eta)=D\eta^2+a_1(c)=D\eta^2+c\lambda+\d-D\lambda^2\ge D\eta^2+pe^{-\lambda cr}=D\eta^2 +k_2(c),
\]
and
\[
k_1(c,\eta)-k_2(c)=D\eta^2 +c\lambda+\d-D\lambda^2-pe^{-\lambda cr}=D\eta^2 +\mu_0, \ \ \mbox{ for } c>c_*,
\]
where $\mu_0:= c\lambda+\d-D\lambda^2-pe^{-\lambda cr}>0$ for $c>c_*$, and
\[
k_1(c_*,\eta)-k_2(c_*)=D\eta^2 +c\lambda_*+\d-D\lambda_*^2-pe^{-\lambda_* c_*r}=D\eta^2, \ \ \mbox{ for } c=c_*.
\]
Then,  from \eqref{12-26-3} we have
\begin{eqnarray}\label{8-2-5}
|e^{-A(\eta)(t+r)}e^{\bar{B}(\eta)t}_r| &\le & e^{-k_1(c,\eta)(t+r)}e^{\bar{k}_2(c,\eta)t}_r \notag \\
&\le& C e^{-\varepsilon_1[k_1(c,\eta)-k_2(c)]t} \notag \\
&=&
\begin{cases}
C e^{-\varepsilon_1(D\eta^2 +\mu_0)t} &\mbox{ for } c>c_*, \\
C e^{-\varepsilon_1 D\eta^2 t} &\mbox{ for } c=c_*.
\end{cases}
\end{eqnarray}
Applying \eqref{8-2-5}, we derive the optimal estimate for $I_1(t)$:
\begin{eqnarray}\label{8-3-4}
I_1(t)&=& \int_{-\infty}^\infty \Big|e^{-A(\eta)(t+r)}e_r^{\bar{B}(\eta)t} \hat{u}^+_0(-r,\eta)\Big|d\eta \notag \\
&\le &
\begin{cases}
C\int_{-\infty}^\infty e^{-\varepsilon_1 D\eta^2 t} e^{-\varepsilon_1\mu_0 t} |{\hat u}_0^+(-r,\eta)| d\eta  \ \ &\mbox{ for } c>c_* \\
C\int_{-\infty}^\infty e^{-\varepsilon_1 D\eta^2 t} |{\hat u}_0^+(-r,\eta)| d\eta  \ \ &\mbox{ for } c=c_*
\end{cases}
\notag \\
&\le &
\begin{cases}
Ce^{-\varepsilon_1\mu_0 t} \|{\hat u}_0^+(-r)\|_{L^\infty}\int_{-\infty}^\infty e^{-\varepsilon_1 D\eta^2 t}  d\eta  \ \ &\mbox{ for } c>c_* \\
C\|{\hat u}_0^+(-r)\|_{L^\infty}\int_{-\infty}^\infty e^{-\varepsilon_1 D\eta^2 t}  d\eta  \ \ &\mbox{ for } c=c_*
\end{cases}
\notag \\
&\le &
\begin{cases}
Ce^{-\mu_1 t} t^{-\frac{1}{2}}\|{ u}_0^+(-r)\|_{L^1}  \ \ &\mbox{ for } c>c_*, \\
Ct^{-\frac{1}{2}}\|{ u}_0^+(-r)\|_{L^1}  \ \ &\mbox{ for } c=c_*,
\end{cases}
\end{eqnarray}
for $\mu_1:=\varepsilon_1 \mu_0 <c\lambda+\d-D\lambda^2-pe^{-c\lambda r}$.

Similarly, we may estimate $I_2(t)$. We note that
\begin{eqnarray}\label{8-3-1}
\sup_{\eta\in R} |A(\eta) {\hat u}^+_0(s,\eta)|&=&\sup_{\eta\in R} \Big|[D\eta^2 + a_1(c)+\mbox{i}a_0(c)\eta] {\hat u}^+_0(s,\eta)\Big| \notag \\
&\le & C\|u^+_0(s)\|_{W^{2,1}(\R)}.
\end{eqnarray}
Thus, we can derive the decay rate for $I_2(t)$ as follows
\begin{eqnarray}\label{8-3-2}
I_2(t)&=&\int^0_{-r} \int_{-\infty}^\infty \Big| e^{-A(\eta)(t-s)}e_r^{\bar{B}(\eta)(t-r-s)}  \Big[\frac{d}{ds}\hat{u}^+_0(s,\eta) + A(\eta)\hat{u}^+_0(s,\eta)\Big] \Big| d\eta ds \notag \\
&\le &  \begin{cases}
C\int^0_{-r}\int_{-\infty}^\infty e^{-\varepsilon_1 D\eta^2 (t-s)} e^{-\varepsilon_1\mu_0 (t-s)} \Big[\frac{d}{ds}\hat{u}^+_0(s,\eta) + A(\eta)\hat{u}^+_0(s,\eta)\Big] \Big| d\eta ds  \\
\hspace{8cm} \mbox{ for } c>c_* \\
C\int^0_{-r}\int_{-\infty}^\infty e^{-\varepsilon_1 D\eta^2 (t-s)} \Big[\frac{d}{ds}\hat{u}^+_0(s,\eta) + A(\eta)\hat{u}^+_0(s,\eta)\Big] \Big| d\eta ds  \\ \hspace{8cm} \mbox{ for } c=c_*
\end{cases}
\notag \\
&\le &
\begin{cases}
Ce^{-\varepsilon_1\mu_1 t}\int^0_{-r} e^{\varepsilon_1\mu_0 s}\sup_{\eta\in R}\Big|\frac{d}{ds}\hat{u}^+_0(s,\eta) + A(\eta)\hat{u}^+_0(s,\eta) \Big|\int_{-\infty}^\infty e^{-\varepsilon_1 D\eta^2 (t-s)}  d\eta  ds \\
\hspace{8cm} \mbox{ for } c>c_* \\
C\int^0_{-r} \sup_{\eta\in R}\Big|\frac{d}{ds}\hat{u}^+_0(s,\eta) + A(\eta)\hat{u}^+_0(s,\eta) \Big|\int_{-\infty}^\infty e^{-\varepsilon_1 D\eta^2 (t-s)}  d\eta  ds  \\
\hspace{8cm}\mbox{ for } c=c_*
\end{cases}
\notag \\
&\le &
\begin{cases}
Ce^{-\varepsilon_1\mu_0 t} \int^0_{-r} (t-s)^{-\frac{1}{2}}[\|{ u_0^+}'(s)\|_{L^1(\R)}+\|u^+_0(s)\|_{W^{2,1}(\R)}
] ds  \\
\hspace{8cm} \mbox{ for } c>c_*, \\
C \int^0_{-r} (t-s)^{-\frac{1}{2}}[\|{ u_0^+}'(s)\|_{L^1(\R)}+\|u^+_0(s)\|_{W^{2,1}(\R)}] ds  \\
\hspace{8cm}  \mbox{ for } c=c_*,
\end{cases}
\notag \\
&\le &
\begin{cases}
Ce^{-\varepsilon_1\mu_1 t} (1+t)^{-\frac{1}{2}} \int^0_{-r} \frac{(1+t)^{\frac{1}{2}}}{(t-s)^{\frac{1}{2}}}[\|{ u_0^+}'(s)\|_{L^1(\R)}+\|u^+_0(s)\|_{W^{2,1}(\R)}] ds  \\
\hspace{8cm}  \mbox{ for } c>c_*, \\
C (1+t)^{-\frac{1}{2}} \int^0_{-r} \frac{(1+t)^{\frac{1}{2}}}{(t-s)^{\frac{1}{2}}}[\|{ u_0^+}'(s)\|_{L^1(\R)}+\|u^+_0(s)\|_{W^{2,1}(\R)}] ds   \\
 \hspace{8cm}  \mbox{ for } c=c_*,
\end{cases}
\notag \\
&\le &
\begin{cases}
Ce^{-\mu_1 t} (1+t)^{-\frac{1}{2}}[\|{ u_0^+}'(s)\|_{L^1(-r,0;L^1(\R))}+\|u^+_0(s)\|_{L^1(-r,0;W^{2,1}(\R))}]   \\
\hspace{8cm}  \mbox{ for } c>c_*, \\
C (1+t)^{-\frac{1}{2}}[\|{ u_0^+}'(s)\|_{L^1(-r,0;L^1(\R))}+\|u^+_0(s)\|_{L^1(-r,0;W^{2,1}(\R))}]  \\
\hspace{8cm}  \mbox{ for } c=c_*.
\end{cases}
\end{eqnarray}
Substituting \eqref{8-3-4} and \eqref{8-3-2} to \eqref{07-29-6}, we obtain the decay rates:
\[
\|u^+(t)\|_{L^\infty(\R)}\le
\begin{cases}
C(1+t)^{-\frac{1}{2}}e^{-\mu_1 t}, & \mbox{ for } c>c_*,\\
C(1+t)^{-\frac{1}{2}}, & \mbox{ for } c=c_*.
\end{cases}
\]
 The proof is complete. $\square$

Let us choose $u_0^+(s,\xi)$ such that   $u_0^+\in L^1([-r,0];W^{2,1}(\R))$ and $\partial_s u^+_0 \in L^1([-r,0];L^1(\R))$ and
\[
u^+_0(s,\xi)\ge  |\tilde{u}_0(s,\xi)|.
\]
Combining Lemmas \ref{lem2.1} and \ref{lem2.2}, we immediately get the convergence rates for $\tilde{u}(t,\xi)$.

\begin{lemma}\label{lem2.3} When $\tilde{u}_0 \in L^1([-r,0];W^{2,1}(\R))$ and $\partial_s \tilde{u}_0 \in L^1([-r,0];L^1(\R))$, then
\begin{equation}
\|\tilde{u}(t)\|_{L^\infty(\R)}\le
\begin{cases}
C(1+t)^{-\frac{1}{2}}e^{-\mu_1 t}, & \mbox{ for } c>c_*,\\
C(1+t)^{-\frac{1}{2}}, & \mbox{ for } c=c_*.
\end{cases}
\label{2.27}
\end{equation}
\end{lemma}

Since $\tilde{u}(t,\xi)=\sqrt{w(\xi)}u(t,\xi)=e^{-\lambda \xi} u(t,\xi)$ for $c\ge c_*$, and $e^{-\lambda \xi}\to 0$ as $\xi\to \infty$,  from \eqref{2.27} we cannot
guarantee any decay estimate for $u(t,\xi)$ at far field $\xi\sim \infty$. In order to get the optimal decay rates for $u(t,\xi)$ in cases of $c>c_*$ and $c=c_*$, let us first investigate
the convergence of $u(t,\xi)$ at far field of $\xi=\infty$.

\begin{lemma}\label{lem2.4} When $e<\frac{p}{\d}\le e^2$ for any $r>0$, or $\frac{p}{\d}>e^2$ but for $r\in (0,\bar{r})$, then there exists a large number $x_0\gg 1$ such that the solution $u(t,\xi)$ of \eqref{2.1} satisfies
\begin{equation}
\sup_{[x_0,\infty)}|u(t,\xi)|\le C e^{-\mu_2 t}, \ \ t>0, \ \ \ c\ge c_*,
\label{2.30}
\end{equation}
for some $0<\mu_2=\mu_2(p,\d,r,b'(v_+))<d$.
\end{lemma}
{\bf Proof}.
Since $u\in C_{unif}(0,\infty)$, namely $\lim_{\xi\to +\infty} u(t,\xi)=u(t,\infty)=:z^+(t)$ exists uniformly for $t\in[-r,\infty)$ and $\lim_{\xi\to +\infty}u_\xi(t,\xi)=\lim_{\xi\to+\infty}u_{\xi\xi}(t,\xi)=0$ are uniformly for $t\in [-r,\infty)$,  let us  take the limits to \eqref{2.1} as $\xi\to\infty$, then we have
\begin{equation}
\begin{cases}
 \dfrac{ d }{dt}z^+(t) + \d z^+(t) -b'(v_+)z^+(t-r)=Q(z^+(t-r)), \\
z^+(s)=z^+_0(s), \ s\in [-r,0],
\end{cases}
\label{2.28}
\end{equation}
where
\[
Q(z^+)=b(v_++z^+)-b(v_+)-b'(v_+)z^+.
\]
As shown in \cite{LLLM}, it is well-known that, when $e<\frac{p}{\d}\le e^2$ for any $r>0$, or $\frac{p}{\d}>e^2$  but with a small time-delay
$0<r<{\bar r}$,  then the above equation \eqref{2.28} satisfies
\begin{equation}
|u(t,\infty)|=|z^+(t)|\le C e^{-\mu_2 t}, \ t>0, c\ge c_*,
\label{2.29}
\end{equation}
for some $0<\mu_2=\mu_2(p,d,r,b'(v_+))<d$, provided with $|z^+_0|\ll 1$.

From the continuity  and the uniform convergence of $u(t,\xi)$ as $\xi \to +\infty$, there exists a large $x_0\gg 1$ such that \eqref{2.29} implies the following convergence immediately
\begin{equation}
\sup_{\xi\in[x_0,+\infty)}|u(t,\xi)|\le C e^{-\mu_2 t}, \ t>0, c\ge c_*,
\label{2.29-1}
\end{equation}
provided $\sup_{\xi\in[x_0,+\infty)}|u_0(s,\xi)|\ll 1$ for $s\in[-r,0]$. Such a smallness for the initial perturbation $u_0$ near $\xi=+\infty$ can be automatically verified, because $\lim_{x\to +\infty} v_0(s,x)=v_+$, which implies $\lim_{\xi\to+\infty} u_0(s,\xi)=\lim_{\xi\to +\infty} [v_0(s,\xi)-\phi(\xi)]
=v_+-v_+=0$ uniformly  for $s\in[-r,0]$. Thus, the proof is complete. $\square$

\begin{lemma}\label{lemma2.8} It holds
\begin{equation}
\sup_{\xi \in (-\infty, x_0]} |u(t,\xi) | \le
\begin{cases}
C(1+t)^{-\frac{1}{2}}e^{-\mu_1 t}, & \mbox{ for } c>c_*,\\
C(1+t)^{-\frac{1}{2}}, & \mbox{ for } c=c_*.
\end{cases}
\label{2.31}
\end{equation}
\end{lemma}
{\bf Proof}. Since
\[
w(\xi)=\begin{cases}
e^{2\lambda|\xi|}, &\mbox{ for } c>c_* \\
e^{2\lambda_*|\xi|}, &\mbox{ for } c=c_*
\end{cases}
\ge
\begin{cases}
e^{2\lambda x_0}, &\mbox{ for } c>c_*, \ \ \xi\le x_0,\\
e^{2\lambda_* x_0}, &\mbox{ for } c=c_*, \ \ \xi\le x_0.
\end{cases}
\]
and $\tilde{u}(t,\xi)=\sqrt{w(\xi)} u(t,\xi)$,
then from \eqref{2.27} we get
\[
\sup_{\xi \in (-\infty, x_0]} |u(t,\xi)| \le
\begin{cases}
C(1+t)^{-\frac{1}{2}}e^{-\mu_1 t}, & \mbox{ for } c>c_*,\\
C(1+t)^{-\frac{1}{2}}, & \mbox{ for } c=c_*.
\end{cases}
\]
The proof is complete. $\square$

{\bf Proof of Proposition \ref{stability-2}}.
Based on Lemmas \ref{lem2.4} and \ref{lemma2.8}, we immediately prove \eqref{alg-stable-2} and \eqref{exp-stable-2} for the convergence rates of $u(t,\xi)$ for $\xi \in \R$, where $0<\mu<\min\{\mu_1,\mu_2\}$. $\square$

  \medskip

{\bf Acknowledgements}
The research of MM was supported in part by
NSERC grant RGPIN 354724-16, FRQNT grant 192571.  The research of KJZ was supported in part by NSFC No.11371082.
The research of QZ was supported in part by NSFC Grant No. 11501514,
Natural Sciences Foundation of Zhejiang Province Grant No. LQ16A010007,
Zhejiang Province Ministry of Education Grand No. Y201533134 and ZSTU Start-Up Fund  No. 11432932611470.


\begin{thebibliography}{99}

\bibitem{Aguerrea-Gomez-Trofimchuk} M. Aguerrea, C. Gomez, and S. Trofimchuk, {\sc On uniqueness of semi-wavefronts,} Math. Ann. 354 (2012), pp. 73--109.

\bibitem{AW1} {\sc D. S. Aronsen and H. F. Weinberger}, {\em Multidimensional nonlinear diffusion arising in
propagation genetics}, Adv. in Math., 30 (1978), pp. 33--76.


\bibitem{Bramson} {\sc M. Bramson}, {\em Convergence of solutions of the Kolmogorov equation to traveling waves}, Mem.
Amer. Math. Soc. 44 (285), (1983).

\bibitem{Chen} {\sc X. Chen}, {\em Existence, uniqueness, and asymptotic stability of traveling waves in nonlocal
evolution equations}, Adv. Differential Equations, 2 (1997), pp. 125--165.

\bibitem{CMYZ} {\sc I-L. Chern, M. Mei, X. Yang and Q. Zhang}, {\em Stability of non-monotone critical traveling waves for reaction-diffusions with time-delay},
J. Differential Equations, 259 (2015), pp. 1503--1541.

\bibitem{E-W} {\sc J.-P. Eckmann and C. E. Wayne}, {\em The nonlinear stability of front solutions for parabolic pde's}, Comm. Math. Phys., 161 (1994), pp. 323--334.

\bibitem{Fang-Zhao} {\sc J. Fang and X.-Q. Zhao}, {\em Existence and uniqueness of traveling waves for non-monotone integral equations with in applications}, J. Differential Equations, 248 (2010), pp. 2199--2226.

\bibitem{Faria-Huang-Wu} {\sc T. Faria, W. Huang and J. Wu}, {\em Traveling waves for delayed reaction-diffusion equations with global response},
 { Proc. Roy. Soc. London}, Ser. A., {462} (2006), pp. 229--261.


 \bibitem{Faria-Trofimchuk1} {\sc T. Faria and S. Trofimchuk}, {\em Nonmonotone traveling waves in single species reaction-diffusion equation with delay},
 { J. Differential Equations},  { 228} (2006), pp. 357--376.

\bibitem{Fife-McLeod} {\sc P. C. Fife and J. B. McLeod}, {\em A phase plane discussion of convergence to travelling fronts
for nonlinear diffusion}, Arch. Rational Mech. Anal., 75 (1980/81), pp. 281--314

 \bibitem{Ga} {\sc T. Gallay}, {\em Local stability of critical fronts in nonlinear parabolic partial differential equations},
Nonlinearity, 7 (1994), pp. 741–764.

\bibitem{G-S} {\sc T. Gallay, A. Scheel},
{\em Diffusive stability of oscillations in reaction-diffusion systems},
Trans. Amer. Math. Soc. 363 (2011), pp. 2571--2598

\bibitem{GBN} {\sc  W. S. C. Gurney, S. P. Blythe and R. M. Nisbet}, {\em  Nicholson's blowflies revisited},  Nature
287 (1980), pp. 17-21.

 \bibitem{GT} {\sc A. Gomez and S. Trofimchuk}, {\em Global continuation of monotone wavefronts}, J. London Math. Soc. 89 (2014), pp 47--68.
 doi:10.1112/jlms/jdt050


\bibitem{GSW} {\sc S.A. Gourley,  J.W.-H. So and J. Wu},
Non-locality of reaction-diffusion equations induced by delay:
biological modeling and nonlinear dynamics, {\it Contemporary
Mathematics, Thematic Surveys}, Kluwer Plenum, Edited by D. V.
Anosov and A. Skubachevskii, (2003) 84-120 (In Russian). Also
appeared in Journal of Mathematical Sciences  {\bf 124} (2004)
5119-5153 (In English).

\bibitem{GW} {\sc S.A. Gourley and J. Wu}, {\em Delayed nonlocal diffusive system in biological invasion and disease spread}, in:
Fields Inst. Commun.,  Vol. { 48}, 2006, pp. 137--200.

\bibitem{Hamel1} {\sc F. Hamel and N. Nadirashvili}, {\em Travelling fronts and entire solutions of the Fisher-KPP
equation in $R^N$ }, Arch. Rational Mech. Anal., 157 (2001), pp. 91--163.

\bibitem{Hamel2}  {\sc F. Hamel and L. Roques}, {\em Uniqueness and stability properties of monostable pulsating fronts},
J. European Math. Soc., 13 (2011), 345--390.

\bibitem{Hou-Li} {\sc  X. Hou and Y. Li}, {\em Local stability of traveling-wave solutions of nonlinear reaction-diffusion
equations}, Discrete Contin. Dyn. Syst., 15 (2006), pp. 681--701.

\bibitem{Huang-Mei-Wang} {\sc R. Huang, M. Mei and Y. Wang}, {\em Planar traveling waves for nonlocal dispersal equation with monostable nonlinearity},
Discrete Contin. Dyn. Syst., 32 (2012), pp. 3621--3649.

\bibitem{HMZZ} {\sc R. Huang, M. Mei, K.-J. Zhang and Q.-F. Zhang}, {\em Stability of oscillating wavefronts for time-delayed nonlocal dispersal equations},
Discrete Contin. Dyn. Syst., 38 (2016), pp. 1331--1353.

\bibitem{KIK} {\sc D. Ya. Khusainov, A. F. Ivanov and I. V. Kovarzh}, {\em Solution of one heat equation with delay},
Nonlinear  Oscillasions, 12 (2009), pp. 260--282.

\bibitem{Ki}{\sc K. Kirchgassner}, {\em On the nonlinear dynamics of travelling fronts}, J. Differential Equations,
96 (1992), pp. 256–-278.

\bibitem{Lau} {\sc K.-S. Lau}, {\em On the nonlinear diffusion equation of Kolmogorov, Petrovsky, and Piscounov},
J. Differential Equations, 59 (1985), pp. 44–-70.

\bibitem{LW} {\sc D. Liang and J. Wu}, {\em Traveling waves and numerical
approximations in a reaction-diffusion equation with nonlocal
delayed effect}, J. Nonlinear Sci., 13 (2003), pp. 289--310.

\bibitem{Liang-Zhao-1} {\sc X. Liang and X.-Q. Zhao}, {\em Asymptotic speeds of spread and traveling waves for monotone
semiflows with applications}, Commun. Pure Appl. Math., 60 (2007), pp. 1–-40. Erratum: 61
(2008), pp. 137–-138.

\bibitem{LLLM} {\sc C.-K. Lin, C.-T. Lin, Y. P. Lin and M. Mei}, {\em Exponential stability of non-monotone traveling waves for Nicholson's blowflies equation},
SIAM J. Math. Anal., 46 (2014), pp. 1053--1084.


\bibitem{Ma} {\sc S. Ma}, {\em Traveling waves for non-local delayed diffusion equations via auxiliary equations},
{ J. Differential Equations}, { 237} (2007), pp. 259--277.

\bibitem{MG} {\sc M. C. Mackey and L. Glass}, {\em Oscillations and chaos in physiological control systems}, Science, 197, Issue 4300, (1977), pp. 287-289.

\bibitem{M-R} {\sc J.F. Mallordy and J.M. Roquejoffre}, {\em A parabolic equation of the KPP type in higher dimensions},
SIAM J. Math. Anal., 26 (1995), pp. 1--20.


\bibitem{Mei-Lin-Lin-So1} {\sc M. Mei, C.-K. Lin, C.-T. Lin and J. W.-H. So},
{\em Traveling wavefronts for time-delayed reaction-diffusion equation:
(I) local nonlinearity},  { J. Differential
Equations}, { 247} (2009), pp. 495--510.


\bibitem{Mei-Lin-Lin-So2} {\sc M. Mei, C.-K. Lin, C.-T. Lin and J. W.-H. So},
{\em Traveling wavefronts for time-delayed reaction-diffusion equation:
(II) nonlocal nonlinearity},  { J. Differential
Equations}, { 247} (2009), pp. 511--529.


\bibitem{Mei-Ou-Zhao} {\sc M. Mei, C. Ou and X.-Q. Zhao}, {\em Global stability of monostable traveling waves for  nonlocal time-delayed reaction-diffusion equations}, { SIAM J.
Math. Anal.}, { 42} (2010), pp. 2762--2790. {\em Erratum}, { SIAM J.
Math. Anal.}, { 44} (2012), pp. 538--540.

\bibitem{Mei-So} {\sc M. Mei and J.W.-H. So},  {\em Stability of strong traveling waves for a nonlocal time-delayed
reaction-diffusion equation}, Proc. Royal Soc. Edinburgh, Section A, 138 (2008), pp. 551–-568.

\bibitem{MSLS} {\sc M. Mei, J.W.-H. So, M.Y. Li and S.S.P. Shen},
{\em Asymptotic stability of traveling waves for the Nicholson's
blowflies equation with diffusion}, { Proc. Royal Soc.
Edinbourgh},  { 134A} (2004), pp. 579--594.



\bibitem{Mei-Wang} {\sc M. Mei and Y. Wang}, {\em Remark on stability of traveling waves for nonlocal Fisher-KPP equations},
{ Intern. J. Num. Anal.  Model.}, Series  B, { 2} (2011), pp. 379-401.

\bibitem{Moet} {\sc H. J. K. Moet}, {\em A note on asymptotic behavior of solutions of the KPP equation}, SIAM J.
Math. Anal., 10 (1979), pp. 728–732.

 \bibitem{Nich1} {\sc A. J. Nicholson}, {\em Competition for food amongst Lucilia Cuprina larvae}, In Proc. VIII International
Congress of Entomology, Stockholm, pp. 277-281 (1948).

\bibitem{Nich2} {\sc A. J. Nicholson}, {\em An outline of the dynamics of animal populations}, Aust. J. Zool. 2 (1954), pp.
9-65.

\bibitem{Nishihara-Wang-Yang} {\sc K. Nishihara, W.-K. Wang and T. Yang}, {\em L$_p$-convergence rate to nonlinear diffusion waves for p-system with damping},
J. Differential Equations, 161 (2000), pp. 191-218.

\bibitem{RTY} {\sc P. Radu, G. Todorova and B. Yordanov}, {\em Decay estimates for wave equations with variable coefficients}, Trans. Amer. Math. Soc. 362 (2010), pp. 2279--2299.

\bibitem{Schaaf} {\sc K. W. Schaaf}, {\em Asymptotic behavior and traveling wave solutions for parabolic functional
differential equations}, Trans. Amer. Math. Soc., 302 (1987), pp. 587–615.

\bibitem{Sattinger} {\sc D. H. Sattinger}, {\em On the stability of waves of nonlinear parabolic systems}, Adv. Math., 22 (1976), pp. 312--355.

\bibitem{S-S} {\sc B. Sandstede and  A. Scheel},
{\em On the stability of periodic travelling waves with large spatial period},
J. Differential Equations. 172 (2001), 134-188


\bibitem{SY} {\sc J. So and Y. Yang}, {\em Direchlet problem for the diffusive Nicholson's blowflies equation},
{ J. Differential Equations}, { 150} (1998), pp. 317--348.


\bibitem{SZ} {\sc J. So and X. Zou}, {\em Traveling waves for the diffusive Nicholson's blowflies equation},
 { Applied Math. Comput.}, {122} (2001), pp. 385--392.

\bibitem{Th-Zh} {\sc H. Thieme and X.-Q. Zhao}, {\em Asymptotic speeds of spread and traveling waves for integral equation
and delayed reaction-diffusion models}, J. Differential Equations, 195 (2003), pp. 430-370.



\bibitem{Todorova} {\sc  G. Todorova and B. Yordanov}, {\em Nonlinear dissipative wave equations with potential}, Contral Methods in PDE-dynamical Systems, Contemp. Math., 426,  Amer. Math. Soc.  (2007), pp. 317--337.

\bibitem{TT} {\sc E. Trofimchuk and S. Trofimchunk}, {\em Admissible wavefront speeds for a single species reaction-diffusion equation with delay}, Discrete Contin. Dyn. Syst. A 20 (2008), pp. 407--423.

 \bibitem{TTT} {\sc E. Trofimchuk, V. Tkachenko and S. Trofimchuk}, {\em Slowly oscillating wave solutions of a single species reaction-diffusion equation with delay}, J. Differential Equations, 245 (2008), pp. 2307--2332.

 \bibitem{Uchiyama} {\sc K. Uchiyama}, {\em The behavior of solutions of some nonlinear diffusion equations for large time},
J. Math. Kyoto Univ., 18 (1978), pp. 453–-508.

\bibitem{VVV} {\sc A. Volpert, Vi. Volpert and Vl. Volpert}, {\em Traveling Wave Solutions of Parabolic Systems},
Transl. Math. Monographs, 140, AMS, Providence, RI, 1994.



 \bibitem{Wang-Yang} {\sc W.-K. Wang and T. Yang}, {\em Existence and stability of planar diffusion waves for 2-D Euler equations with damping},
J. Differential Equations, 242 (2007), pp. 40--71.

\bibitem{Xin} {\sc J. Xin}, {\em Front propagation in heterogeneous media}, SIAM Review, 42 (2000), pp. 161--230.








\end{thebibliography}
 \end{document}